\documentclass[11pt,leqno]{amsart}
\usepackage{amssymb}
\usepackage{mathrsfs}
\usepackage{color}
\usepackage{soul}
\usepackage{dsfont}
\usepackage{fancyhdr}
\numberwithin{equation}{section}

\def\sgn{\mathop{\rm sgn}}

\newcommand\R{{\mathbb R}}

\newcommand\C{{\mathbb C}}

\newcommand\PPP{{\bf P}}

\def\AA{{\mathcal A}}
\def\BB{{\mathcal B}}

\def\DD{{\mathcal D}}

\def\JJ{{\mathcal J}}

\def\MM{{\mathcal M}}
\def\NN{{\mathcal N}}
\def\OO{{\mathcal O}}
\def\PP{{\mathcal P}}
\def\QQ{{\mathcal Q}}
\def\RR{{\mathcal R}}

\def\XX{{\mathcal X}}

\def\ZZ{{\mathcal Z}}

\def\PPP{{\mathbf P}}

\def\BBB{{\mathscr{B}}}
\def\III{{\mathscr{I}}}
\def\LLL{{\mathscr{L}}}

\def\Ppp{{\mathbf P}}

\def\eps{{\varepsilon}}

\newtheorem{theo}{Theorem}[section]

\newtheorem{lem}[theo]{Lemma}

\newcommand{\beqn}{\begin{equation}}
\newcommand{\eeqn}{\end{equation}}
\newcommand{\bear}{\begin{eqnarray}}
\newcommand{\eear}{\end{eqnarray}}
\newcommand{\bean}{\begin{eqnarray*}}
\newcommand{\eean}{\end{eqnarray*}}
\newcommand{\bal}{\begin{aligned}}
\newcommand{\eal}{\end{aligned}}

\newcommand{\be}{\begin{equation}}
\newcommand{\ee}{\end{equation}}
\newcommand{\ba}{\begin{aligned}}
\newcommand{\ea}{\end{aligned}}

\newcommand{\ds}{\displaystyle}

\def\signqw{\bigskip \begin{center} {\sc Qilong Weng\par\vspace{3mm}
Universit\'e Paris-Dauphine \par
CEREMADE, UMR CNRS 7534\par
Place du Mar\'echal de Lattre de Tassigny
75775 Paris Cedex 16\par
FRANCE\par\vspace{3mm}
e-mail:} \tt{weng@ceremade.dauphine.fr} \end{center}}

\begin{document}

\title{General time elapsed neuron network model: well-posedness
and strong connectivity regime}

\author{Q. WENG}

\begin{abstract}
For large fully connected neuron networks, we study the dynamics of
homogenous assemblies of interacting neurons described by time
elapsed models, indicating how the time elapsed since the last
discharge construct the probability density of neurons. Through the
spectral analysis theory for semigroups in Banach spaces developed
recently in \cite{Mbook*,MS}, on the one hand, we prove the existence 
and the uniqueness of the weak solution in the whole connectivity
regime as well as the parallel results on the long time behavior of 
solutions obtained in \cite{MSWQ} under general assumptions on the 
firing rate and the delay distribution. On the other hand, we extend
those similar results obtained in \cite{PK1,PK2} in the case without
delay to the case taking delay into account and both in the weak and
the strong connectivity regime with a particular step function firing
rate. 
\end{abstract}

\maketitle

\begin{center} {\bf Version of \today}
\end{center}

\vspace{0.3cm}

\bigskip
\bigskip


\tableofcontents


\smallskip
\noindent \textbf{Keywords.} Neuron networks, time elapsed dynamics,
semigroup spectral analysis, strong connectivity, exponential
asymptotic stability.

\maketitle

\section{Introduction}
\label{sec:Intro}
The information transmission and processing mechanism in the nervous
systems relies on the quantity of electrical pulses as the reflect to
incoming stimulations, during which the neurons experience a period of
recalcitrance called discharge time before reactive. In this work, we
shall focus on the model describing the neuronal dynamics in
accordance with this kind of discharge time which has been introduced
and studied in \cite{GKbook,PK1,PK2}. In order to show the response to the
recovery of the neuronal membranes after each discharge, the model 
consider an instantaneous firing rate depending on the time elapsed
since last discharge as well as the inputs of neurons. This sort of
models are also regarded as a mean field limit of finite number of
neuron network models referred to \cite{MR3311484,FL*,RT*,Q*}.

For a local time (or internal clock) $x\ge0$ corresponding to the
elapsed time since the last discharge, we consider the dynamic of the
neuronal network with the density number of neurons $f=f(t,x)\geq0$ 
in state $x \ge 0$ at time $t \ge 0$, given by the following nonlinear 
time elapsed (or of age structured type) evolution equation
\begin{subequations}\label{eq:ASM}
 \begin{align}
   &\partial_t f=-\partial_x f-a(x,\eps \,  m(t))f=:\mathcal{L}_{\eps m(t)}f,\\
   &f(t,0)=p(t), \ \ f(0,x)=f_0(x),
 \end{align}
\end{subequations}
where $a(x,\eps \, \mu)\ge0$ denotes the firing rate of a neuron in
state $x$ and in an environment $\mu\ge0$ formed by the global 
neuronal activity with a network connectivity parameter $\eps\ge0$
corresponding to the strength of the interactions. The total density 
of neurons $p(t)$ undergoing a discharge at time $t$ is defined through
$$
p(t):=\mathcal{P}[f(t); m(t)],
$$
where
$$
\PP [g,\mu] = \PP_\eps[g,\mu] :=  \int_0^\infty a(x,\eps \mu) g(x)\mathrm{d}x,
$$
while the global neuronal activity $m(t)$ at time $t\geq0$ taking into
account the interactions among the neurons resulting from earlier
discharges is given by
\beqn\label{eq:delay}
m(t):=\int_0^\infty p(t-y)b(\mathrm{d}y),
\eeqn
where the delay distribution $b$ is a probability measure considering
the persistence of the electric activity to those discharges in the
network. In the sequel, we will consider the two following situations
respectively:

$\bullet$ The {\it case without delay},  when $b=\delta_0$ then $m(t) = p(t)$.

$\bullet$ The {\it case with delay},  when $b$ is smooth.

\smallskip
Observe that the solution $f$ of the time elapsed equation
\eqref{eq:ASM} satisfies
$$
\frac{\mathrm{d}}{\mathrm{d}t}\int_0^{\infty}f(t,x)\mathrm{d}x=f(t,0)-
\int_0^{\infty}a(x,\eps m(t))f(t,x)\mathrm{d}x=0,
$$
in both the cases, which implies the conservation of the total density
number of neurons (also called {\it mass} in the sequel) permitting us
to normalize it to be $1$. Then we assume in the sequel 
$$
\langle f(t,\cdot) \rangle = \langle f_0 \rangle = 1,
\quad \forall t\geq0, \quad \langle g \rangle := \int_0^{\infty} g(x)\mathrm{d}x.
$$

\smallskip
We define a couple $(F_\eps,M_\eps)$ as a corresponding steady state,
which satisfy
\begin{subequations}\label{eq:StSt}
 \begin{align} \label{eq:StSt1}
   &0=-\partial_x F_\eps - a(x,\eps \,  M_\eps) F_\eps= \mathcal{L}_{\eps M_\eps} F_\eps,\\
   &F_\eps(0)= M_\eps, \quad \langle F_\eps \rangle = 1.
 \end{align}
\end{subequations}
Noticing that the associated network activity and the discharge
activity are equal constants for a steady state as $\langle b\rangle=1$.

\smallskip
Our main purpose in this paper is to prove the
existence and uniqueness of the solution to the time elapsed
evolution equation \eqref{eq:ASM} no matter which $\eps>0$.
Furthermore, we obtain the exponential asymptotic stability in strong
connectivity regime, which is a range of connectivity parameter 
$\eps\in[\eps_1,\infty)$, with $\eps_1$ large enough, such that the
equations \eqref{eq:ASM} and \eqref{eq:StSt} do not possess intense
nonlinearity. We are also able to extend the result in \cite{PK1,PK2} in 
the case without delay for a {\it step function} firing rate $a$ to the 
case with delay rather than the indescribable stability. In order to 
conclude those results, it is necessary to give the following 
mathematical assumptions on the firing rate $a$  and on the 
delay distribution $b$.

\smallskip
We make the physically reasonable assumptions
\beqn\label{hyp:a1}
\partial_x a\geq0, \ \ a'=\partial_\mu a\geq0,
\eeqn
\beqn\label{hyp:a2}
0<a_0:=\lim_{x\rightarrow\infty}a(x,0)\leq\lim_{x,\,\mu\rightarrow\infty}a(x,\mu)=:a_1<\infty,
\eeqn
one particular example of the firing rate is the "step function" one
\beqn\label{hyp:s1}
a(x,\mu) = {\bf 1}_{x > \sigma(\mu)},\quad \sigma' \le 0,
\eeqn
\beqn\label{hyp:s2}
\sigma(0)=\sigma_+, \quad\sigma(\infty)=\sigma_-<\sigma_+<1,
\eeqn 
associated with some continuity assumption
\beqn\label{hyp:a3}
a \in W^{1,\infty}(\R^2_+),
\eeqn
or particularly
\beqn\label{hyp:s3}
\quad \sigma,\sigma^{-1}\in W^{1,\infty}(\R_+), 
\eeqn
In the strong connectivity regime, we consider the decay assumption on the two cases, for a.e. $x\ge0$,
\beqn\label{hyp:a4}
\eps\sup_{x\ge0}\partial_\mu a(x,\eps\mu)\to0, \quad\mathrm{as} \ \eps\to\infty,
\eeqn
\beqn\label{hyp:s4}
\limsup\limits_{\eps\to\infty}\sup\limits_{\mu\in[1-\sigma_+,1]}\eps|\sigma'(\eps\mu)|=0.
\eeqn
In the case with delay, we assume that the delay distribution
$b(\mathrm{d}y) = b(y)\mathrm{d}y$ has the exponential bound and
satisfies the smoothness condition
\begin{equation}\label{hyp:del}
\exists\delta>0, \quad\int_0^\infty e^{\delta
y}\, (b(y) +  |b'(y)|) \, \mathrm{d}y <\infty.
\end{equation}

The above assumptions permit the existence and uniqueness of the 
solution to the nonlinear problem \eqref{eq:ASM} thanks to the Banach 
fixed-point theorem.

\begin{theo}\label{th:EaU}
On the one hand, for a smooth firing rate, assuming \eqref{hyp:a1}-\eqref{hyp:a2}-\eqref{hyp:a3}-\eqref{hyp:a4}, then for any $f_0\in L^1(\R_+)$ and any $\eps>0$, there exists a unique nonnegative and mass conserving weak solution $f\in C(\R_+;L^1(\R_+))$ to the evolution
equation \eqref{eq:ASM} for some functions $m,\,p\in C([0,\,\infty))$.
On the other hand, for a step function firing rate, assuming \eqref{hyp:s1}-\eqref{hyp:s2}-\eqref{hyp:s3}-\eqref{hyp:s4}, then for any $f_0\in L^1\cap L^\infty(\R_+)$ and $0\le f_0\le1$, the parallel conclusion holds with the unique solution $(f,m)$ satisfying that
$$
0\leq f(t,x)\leq1, \quad\forall t,x\geq0,
$$
$$
1-\sigma_+\le m(t)\le 1, \quad\forall t\geq0.
$$
\end{theo}
For any $\eps>0$, there also exits a corresponding steady state, 
which is unique additionally in the weak or strong connectivity regime.

\begin{theo}\label{th:SS}
On the one hand, for a smooth firing rate, assuming \eqref{hyp:a1}-\eqref{hyp:a2}-\eqref{hyp:a3}-\eqref{hyp:a4}, then for any $\eps\ge0$, there exists at least one pair of solutions $(F_\eps(x),M_\eps)\in W^{1,\infty}(\R_+)\times\R_+$ to the stationary problem \eqref{eq:StSt} such that
\beqn\label{ineq:EoF}
0\le F_\eps(x)\lesssim e^{{-a_0\over2}x},\quad |F'_\eps(x)|\lesssim e^{{-a_0\over2}x},\quad x\ge0.
\eeqn
On the other hand, for a step function firing rate, the existence of the steady state holds under the assumptions \eqref{hyp:s1}-\eqref{hyp:s2}-\eqref{hyp:s3}-\eqref{hyp:s4}, which satisfy
$$
0\le F_\eps\le1,\quad 1-\sigma_+\le M_\eps\le1.
$$
Moreover, there exists $\eps_0>0$ small enough or $\eps_1 > 0$ large enough, such that the above steady state is unique for any $\eps\in[0,\eps_0)\cup(\eps_1,+\infty]$.
\end{theo}

\smallskip
We conclude the exponential long time stability in the strong connectivity regime combined with those in the weak one in \cite{MSWQ} as our main result.

\begin{theo}\label{th:MR} We assume that the firing rate $a$
 satisfies \eqref{hyp:a1}-\eqref{hyp:a2}-\eqref{hyp:a3}-\eqref{hyp:a4} or \eqref{hyp:s1}-\eqref{hyp:s2}-\eqref{hyp:s3}-\eqref{hyp:s4}. We also assume that the  delay distribution  $b$  satisfies $b = \delta_0$ or  \eqref{hyp:del}. There exists $\eps_0>0 \ (\eps_1> 0)$, small (large) enough, such that for any
$\eps\in(0, \eps_0) \ (\eps\in(\eps_1,+\infty))$ the steady state $(F_\eps,M_\eps)$ is unique. There also exist some constants $\alpha< 0$, $C \ge 1$ and  $\eta > 0$ (besides $\zeta_\eps\to0$ as $\eps\to\infty$) such that for any connectivity parameter $\eps\in(0,\eps_0) \ (\eps\in(\eps_1,+\infty))$  and for any unit mass initial datum $0 \le f_0\in L^1$ or in the case of step function firing rate for any unit mass initial datum $f_0\in L^p$, $1\le p\le\infty$, satisfying $0\le f_0\le1$ additionally, such that $\| f_0 - F_\eps \|_{L^1} \le \eta/\eps \ (\le\eta/\zeta_\eps)$, then the (unique, positive and mass conserving) solution $f$ to the evolution equation \eqref{eq:ASM} satisfies
$$
\|f(t,.)-F_\eps \|_{L^1}\leq C  e^{\alpha t}, \qquad \forall \, t \ge 0 .
$$
\end{theo}

\smallskip
In order to study the asymptotic convergence to an equilibrium for the
homogeneous inelastic Boltzman equation, the strategy of
``perturbation of semigroup spectral gap'' is first introduced in
\cite{MMcmp}. Inspired by its recent application to a neuron network equation in
\cite{MQT*}, we linearize the equation around a stationary state
$(F_\eps,M_\eps,M_\eps)$ on the variation $(g,n,q)=(f,m,p)-(F_\eps,M_\eps,M_\eps)$, 
such that
\begin{subequations}\label{eq:ASMlin1}
 \begin{align}
   &\partial_t g=-\partial_x g- a(x,\eps \, M_\eps) g - n(t) \, \eps (\partial_\mu a)(x,\eps M_\eps) F_\eps,  \\
   &g(t,0)=q(t), \ \ g(0,x)= g_0(x),
    \end{align}
\end{subequations}
with
\beqn\label{eq:ASMlin2}
q(t) =   \int_0^\infty a(x,\eps \, M_\eps) g \,\mathrm{d}x +  n(t) \, \eps \int_0^\infty  (\partial_\mu a)(x,\eps M_\eps) F_\eps \,  \mathrm{d}x
\eeqn
and
\beqn\label{eq:ASMlin3}
  n(t):=\int_0^\infty q(t-y)b(\mathrm{d}y),
\eeqn
while for the step function firing rate \eqref{hyp:s1} in \cite{PK1,PK2}, it is impossible to linearize the model because of the failure in meeting the condition \eqref{hyp:a3}. Insteadly, we introduce another more concise linear equation around the steady state on the variation $(g,n,p)$, which writes
\begin{subequations}\label{eq:ASMFRlin1}
 \begin{align}
   &\partial_t g=-\partial_x
   g-g\mathbf{1}_{x>\sigma_\eps},  \\
   &g(t,0)=q(t), \quad g(0,x)= g_0(x),
    \end{align}
\end{subequations}
where here and below we note $\sigma_\eps:=\sigma(\eps M_\eps)$ for
simplicity, with
\beqn\label{eq:ASMFRlin2}
q(t) = \PP[g,M_\eps]= \int_0^\infty g\mathbf{1}_{x>\sigma_\eps}\mathrm{d}x.
\eeqn
By regarding the boundary condition as a source term, we construct the linear generator $\Lambda_\eps$ and the associated semigroup $S_{\Lambda_\eps}$ respectively from the above linear equations to apply the spectral analysis. As in \cite{Voigt80,GMM*,MS,Mbook*}, we split the operator $\Lambda_\eps$ into two parts, one of which is $\alpha$-hypodissipative, $\alpha<0$, denoted by $\BB_\eps$ while the other one is bounded and $\BB_\eps$- power regular, denoted as $\AA_\eps$. Benefiting from this split, the semigroup $S_{\Lambda_\eps}$ admits a finite dimensional dominant part, thanks to a particular version of the Spectral Mapping Theorem in \cite{MS,Mbook*} and the Weyl's Theorem in \cite{Voigt80,GMM*,MS,Mbook*}. As $\eps\to\infty$, the limited semigroup $S_\infty$ becomes positive because of the vanishment of the items with $n(t)$ in \eqref{eq:ASMlin1} and \eqref{eq:ASMlin2}, which permits the Krein-Rutman Theorem established in \cite{MS,Mbook*} to imply that the steady state $(F_\infty,M_\infty)$ possesses the exponential stability. And so does the stationary state $(F_\eps,M_\eps)$ in the strong connectivity regime after a perturbative argument developed in \cite{MMcmp,T*,Mbook*}. Then we extend the exponential stability to our main result Theorem~\ref{th:MR} in the case without delay by the analysis on the rest term of the linear equation \eqref{eq:ASMlin1} or \eqref{eq:ASMlin1} compared to the original nonlinear equation. As for the delay case, we replace the delay equation \eqref{eq:delay} by a simple age equation to form an autonomous system with the linear equation \eqref{eq:ASMlin1} or \eqref{eq:ASMlin1} to generate a semigroup and follow the same strategy.

Actually, the previous works \cite{PK1, PK2} shows the asymptotic stability with simpler explicitly expression and appropriate norm benefiting from the choice of the step function firing rate \eqref{hyp:s1}. Different from those approach, we take consideration of more realistic and flexible firing rate with more abstract method, which allows us to obtain the dissipativity of the corresponding linear operators without the explicitly exhibited norm. In particular, we are able to establish the existence and uniqueness of weak solutions first as the complement of the results in \cite{MSWQ} in the weak connectivity regime and adapt our approach to the strong connectivity regime as well as to the step function firing rate to generalize the stability results obtained in \cite{PK1,PK2} in the case without delay to the case considering of the delay term.

\medskip
This article is organized by the following plan. In Section~\ref{sec:ExitSteS}, we demonstrate the existence and uniqueness of the solution and the stationary state result. In the strong connectivity regime, we introduce the strategy and establish Theorem~\ref{th:MR} in the case without delay in Section~\ref{sec:WithoutDelay}, meanwhile the case with delay in section~\ref{sec:WithDelay}. In section~\ref{sec:StepFctFR}, we establish Theorem~\ref{th:MR} again both in the weak regime and the strong regime for the model with the step function firing rate \eqref{hyp:s1}.

\section{Existence and the steady state}\label{sec:ExitSteS}

\subsection{Existence of the solution}

To establish the existence of a solution to \eqref{eq:ASM}, we are going to
apply a fixed point argument with the benefit of the following lemma.

\begin{lem}\label{lem:ext}
Assuming the firing rate $a$ satisfies \eqref{hyp:a1}-\eqref{hyp:a2}-\eqref{hyp:a3}-\eqref{hyp:a4} for any $m\in L^\infty([0,\,T])$ and $f\in C(\R_+;L^1(\R_+))$ satisfying the
equation \eqref{eq:ASM}, consider the application $\JJ$:
$L^\infty([0,\,T])\rightarrow L^\infty([0,\,T])$,
$$
\JJ(m)(t):=\int_0^t p(t-y)b(\mathrm{d}y), \quad with\quad
p(t)=\int_0^\infty a(x,\eps\,m(t))f(x,t)\mathrm{d}x,
$$
then there exist $T>0$ and $0<C<1$ such that the estimate
\beqn\label{eq:EoJ}
\|\JJ(m_1)-\JJ(m_2)\|_{L^\infty([0,\,T])}\leq\,C\|m_1-m_2\|_{L^\infty([0,\,T])}
\eeqn
holds for all $(m_1,\,m_2)\in\,L^\infty([0,\,T])$ and for any $\eps>0$.
\end{lem}
\proof[Proof of Lemma \ref{lem:ext}] Thanks to the method of characteristics
with respect to $x$ or $t$, a solution $f(t,x)$ to the equation
\eqref{eq:ASM} can be expressed as
\beqn\label{eq:chf1}
f(x,t)=f_0(x-t) e^{-\int_0^t
  a(s+x-t,\eps\,m(s))\mathrm{d}s}, \qquad \forall\, x\geq t
\eeqn
or
\beqn\label{eq:chf2}
f(x,t)=p(t-x) e^{-\int_0^x
  a(s,\eps\,m(s+t-x))\mathrm{d}s}, \qquad \forall\, x\leq t.
\eeqn
We denote $f_i$, $i=1,\,2$, two solutions to the equation
\beqn
\begin{cases}
\partial_t f_i(x,t)+\partial_x f_i(x,t)+a(x,\eps\,m_i(t))f_i(x,t)=0,\\
f_i(0,t)=p_i(t)=\ds\int_0^\infty a(x,\eps\,m_i(t))f_i(x,t)\mathrm{d}x,
\end{cases}
\eeqn
with the same initial data $f_0$. We rewrite $p_i(t)$ corresponding to two expressions \eqref{eq:chf1} and \eqref{eq:chf2} as two global activity functions
\bean
p_i(t) &=& \int_0^t a(x,\eps\,m_i(t)) p_i(t-x) e^{-\int_0^x
  a(s,\eps\,m_i(s+t-x))\mathrm{d}s}\mathrm{d}x\\
&+&\int_t^\infty a(x,\eps\,m_i(t))f_0(x-t) e^{-\int_0^t
  a(s+x-t,\eps\,m_i(s))\mathrm{d}s}\mathrm{d}x.
\eean
We split $p_1(t)-p_2(t)$ into two items $I_1(t)$ and $I_2(t)$
with
\bean
I_1(t) &=& \int_0^t\Big(p_1(t-x)a(x,\eps\,m_1(t))e^{-\int_0^x
  a(s,\eps\,m_1(s+t-x))\mathrm{d}s}\\
 && -p_2(t-x)a(x,m_2(t))e^{-\int_0^x a(s,\eps\,m_2(s+t-x))\mathrm{d}s}\Big)\mathrm{d}x
\eean
and $I_2(t)$ as the remainder.
In order to control the first item, we divide it into three parts as
$I_1(t)=I_{1,1}(t)+I_{1,2}(t)+I_{1,3}(t)$, where
\bean
&&I_{1,1}(t)=\int_0^t(p_1-p_2)(t-x)a(x,\eps\,m_1(t))e^{-\int_0^x
  a(s,\eps\,m_1(s+t-x))\mathrm{d}s}\mathrm{d}x,\\
&&I_{1,2}(t)=\int_0^t
p_2(t-x)\big(a(x,\eps\,m_1(t))-a(x,\eps\,m_2(t))\big) e^{-\int_0^x
  a(s,\eps\,m_1(s+t-x))\mathrm{d}s}\mathrm{d}x,\\
&&I_{1,3}(t)=\int_0^t p_2(t-x)a(x,m_2(t))\big(e^{-\int_0^x
  a(s,\eps\,m_1(s+t-x))\mathrm{d}s}-e^{-\int_0^x a(s,\eps\,m_2(s+t-x))\mathrm{d}s}\big)\mathrm{d}x.
\eean
Clearly, we have the estimates
$$
\|I_{1,1}\|_{L^\infty([0,T])}\leq a_1T\|p_1-p_2\|_{L^\infty([0,T])}
$$
and
$$
\|I_{1,2}\|_{L^\infty([0,T])}\leq \eps a_1\|\partial_\mu a\|_{L_x^\infty}T \|m_1-m_2\|_{L^\infty([0,T])}.
$$
Since there exists a constant $C$ such that
\bean
&&\Big|e^{-\int_0^x
  a(s,\eps\,m_1(s+t-x))\mathrm{d}s}-e^{-\int_0^x
  a(s,\eps\,m_2(s+t-x))\mathrm{d}s}\Big|\\
&\leq&  C\int_0^x|a(s,\eps\,m_1(s+t-x))-a(s,\eps\,m_2(s+t-x))|\mathrm{d}s,
\eean
which leads to the estimate
$$
\|I_{1,3}\|_{L^\infty([0,T])}\leq \eps a_1^2\|\partial_\mu a\|_{L_x^\infty}\frac{T^2}{2}\|m_1-m_2\|_{L^\infty{[0,T]}}.
$$
From the assumption \eqref{hyp:a4}, there exists $\eps_1$ large enough such that
$\eps\|a'\|_{L^\infty_x}\leq 1$, for any $\eps\in[\eps_1,+\infty)$. Denoting $\eta:=\max\{\eps_1,\,1\}$, we deduce that
\beqn\label{ineq:EI1}
\|I_1\|_{L^\infty([0,T])}\leq a_1 T \|p_1-p_2\|_{L^\infty([0,T])}+\eta(C_1
T^2+C_2 T) \|m_1-m_2\|_{L^\infty([0,T])},
\eeqn
for any $\eps>0$. On the other hand, we have the item $I_2(t)$ as
\bean
I_2(t) &=& \int_t^\infty f_0(x-t)\Big(a(x,\eps\,m_1(t))e^{-\int_0^t
  a(s+x-t,\eps\,m_1(s))\mathrm{d}s}\\
 & & -a(x,\eps\,m_2(t))e^{-\int_0^t a(s+x-t,\eps\,m_2(s))\mathrm{d}s}\Big)\mathrm{d}x\\
 &=& \int_0^\infty f_0(x)\Big(a(x+t,\eps\,m_1(t))e^{-\int_0^t
  a(s+x,\eps\,m_1(s))\mathrm{d}s}\\
 & & -a(x+t,\eps\,m_2(t))e^{-\int_0^t a(s+x,\eps\,m_2(s))\mathrm{d}s}\Big)\mathrm{d}x
\eean
Clearly, we have
$$
a(x+t,\eps\,m(t))=\frac{\mathrm{d}}{\mathrm{d}t}\int_0^t a(s+x,\eps\,m(s))\mathrm{d}s,
$$
which implies
\bean
 && \Big|a(x+t,\eps\,m_1(t))e^{-\int_0^t a(s+x,\eps\,m_1(s))\mathrm{d}s}
 -a(x+t,\eps\,m_2(t))e^{-\int_0^t a(s+x,\eps\,m_2(s))\mathrm{d}s}\Big|\\
 &\le& C\int_0^t |a(s+x,\eps\,m_1(s))-a(s+x,\eps\,m_2(s))|\mathrm{d}s,
\eean
for some constant $C>0$. Similarly to estimating $I_1$, we deduce that
\beqn\label{ineq:EI2}
\|I_2\|_{L^\infty([0,T])}\leq C_3\eta T \|m_1-m_2\|_{L^\infty([0,T])}.
\eeqn
From the above estimates \eqref{ineq:EI1} and \eqref{ineq:EI2}, it turns out that
$$
\|p_1-p_2\|_{L^\infty([0,T])} \leq a_1T
\|p_1-p_2\|_{L^\infty([0,T])}+\eta T(C_1 T+C_2') \|m_1-m_2\|_{L^\infty([0,T])},
$$
which implies
\beqn\label{ineq:ep}
\|p_1-p_2\|_{L^\infty([0,T])}\leq \eta C T \|m_1-m_2\|_{L^\infty([0,T])}
\eeqn
when $a_1T$ less than $1$. Hence, in the case without delay, we have
$$
\|\JJ(m_1)-\JJ(m_2)\|_{L^\infty([0,T])}=\|p_1-p_2\|_{L^\infty([0,T])}\leq \eta C T \|m_1-m_2\|_{L^\infty([0,T])},
$$
while taking the delay into account with the fact that
$$
\JJ(m_1)(t)-\JJ(m_2)(t)=\int_0^t(p_1-p_2)(t-y)b(\mathrm{d}y),
$$
we obviously deduce
$$
\|\JJ(m_1)-\JJ(m_2)\|_{L^\infty([0,T])}\leq \eta CT^2\|m_1-m_2\|_{L^\infty([0,T])}
$$
from \eqref{ineq:ep}. By taking $T$ small enough such that $\eta CT^2<\eta C T<1$, we finally attain our estimate \eqref{eq:EoJ}.
\endproof

\proof[Proof of Theorem \ref{th:EaU}] From Lemma
\ref{lem:ext}, for any $\eps>0$, there is a $T>0$ which does not depend
upon the initial data such that the application $\JJ$ admits a unique 
fixed point $m(t)$ on $[0,T]$ then the corresponding $f(t,x)$ on 
$[0,T]\times\R^d$, which is the unique solution to the equation 
\eqref{eq:ASM}, according to the Banach-Picard fixed point theorem.
Iterating on $T$, we deduce the global existence and uniqueness 
of the solution $(f,m)$ to equation \eqref{eq:ASM}.
\endproof

\subsection{The stationary problem}

Now we present the proof of the steady state in the strong
connectivity regime.
\proof[Proof of Theorem \ref{th:SS}]  {\sl Step 1. Existence.}
From the assumption \eqref{hyp:a2}, we deduce that for any $x\ge0$,
$\mu\ge0$, there exists $x_0\in[0,\infty)$ such that
$a(x,\mu)\geq\frac{a_0}{2}$. Denoting
$$
A(x,\mu):=\int_0^x a(y,\mu)\mathrm{d}y, \quad \forall ±, x,m \ge 0,
$$
we naturally estimate that
\begin{equation}\label{ineq:EoA}
    \frac{a_0}{2}(x-x_0)_+\leq A(x,\mu)\leq a_1 x, \quad \forall \, x\geq0, \ \mu \ge 0.
\end{equation}
For any $m\ge0$, the equation \eqref{eq:StSt1} can be solved by
$$
F_{\eps,m}(x):=T_m e^{-A(x,\eps m)},
$$
whose mass conservation gives
$$
T_m^{-1}=\int_0^{\infty}e^{-A(x,\eps m)}\mathrm{d}x.
$$
Then the existence of the solution is equivalent to find $m = M_\eps$ satisfying
$m = F_{\eps,m}(0) = T_m$. Considering 
$$
 \Psi(\eps,m)=m
T_{m}^{-1} := m\int_0^{\infty}e^{-A(x,\eps m)}\mathrm{d}x,
$$
it is merely necessary to find $M_\eps\ge0$ such that
\begin{equation}\label{eq:mnc}
\Psi(\eps,M_\eps)=1.
\end{equation}
From the Lebesgue dominated convergence theorem, the function $\Psi(\eps,.)$ is
continuous. In addition to the fact that $\Psi(0)=0$ and $\Psi(\infty)=\infty$, 
the intermediate value theorem implies the existence immediately. The inequality \eqref{ineq:EoA} shows the estimates \eqref{ineq:EoF} clearly.

\smallskip\noindent   {\sl Step 2. Uniqueness in the strong connectivity regime.}
Obviously,
$$
M_\infty:=(\int_0^{\infty}e^{-A(x,\infty)}\mathrm{d}x)^{-1}\in(0,\infty)
$$ 
is the unique solution to $\Psi(\infty,M_\infty)=1$. It is clear that
$$
\frac{\partial}{\partial
m}\Psi(\eps,m)=\int_0^{\infty}e^{-A(x,\eps
m)} \big(1 -m\int_0^x\eps\partial_\mu a(y,\eps m)\mathrm{d}y\big)\mathrm{d}x,
$$
is continuous with respect to the two variables because of
\eqref{hyp:a3}, which implies that $\Psi\in C^1$. Coupled with that
$$
\frac{\partial}{\partial
m}\Psi(\eps,m)|_{\eps=\infty}=\int_0^{\infty}e^{-A(x,\infty)}\mathrm{d}x>0,
$$
we conclude from the implicit function theorem that there exists 
$\eps_1>0$, large enough, such that the equation \eqref{eq:mnc} 
has a unique solution for any $\eps\in(\eps_1,+\infty]$.
\endproof

\section{Case without delay}
\label{sec:WithoutDelay}

In this section, we conclude our main result Theorem~\ref{th:MR}
gradually in the case without delay.

\subsection{Linearized equation and structure of the spectrum}
We introduce the linearized equation on the variation $(g,n)$ around 
the steady state $(F_\eps,M_\eps)$ given by
\bean
   &&\partial_t g+\partial_x g+a_\eps g+a'_\eps F_\eps n(t)=0,\\
   &&g(t,0)=n(t)=\int_0^{\infty}(a_\eps g+a'_\eps F_\eps n(t))\, \mathrm{d}x,
   \quad g(0,x)=g_0(x),
\eean
with notes $a_\eps :=a(x,\eps M_\eps)$ and $a'_\eps := \eps \, (\partial_\mu
a)(x,\eps M_\eps)$ for simplification. According to the assumption
\eqref{hyp:a4}, there exists $\eps_1 > 0$, large enough, such that
$$
\forall \, \eps \in (\eps_1,\infty) \qquad \kappa := \int_0^{\infty} a'_\eps F_\eps \mathrm{d}x < 1,
$$
permitting to define
\beqn\label{def:Nepsg}
n(t)=\MM_\eps[g] := (1-\kappa)^{-1} \int_0^\infty a_\eps g\, \mathrm{d}x.
\eeqn
We consider the operator $L_\eps$ to the above linearized euqation given by
$$
L_\eps g:=\partial_x g-a_\eps g-a'_\eps F_\eps\MM_\eps[g]
$$
in domain
$$
D(L_\eps):=\{g\in W^{1,1}(\R_+), g(0)=\MM_\eps[g]\}
$$
generating the semigroup $S_{L_\eps}$ on space $X:=L^1(\R_+)$. Then for any initial datum $g_0\in X$, the weak solution of the linearized equation is given by $g(t)=S_{L_\eps}(t)g_0$.
By regarding the boundary condition as a source term, we rewrite the linearized equation as
\beqn\label{eq:ASMDL}
\partial_t g=\Lambda_\eps g :=-\partial_x g-a_\eps g-a'_\eps F_\eps\MM_\eps[g]+\delta_{x=0}\MM_\eps[g],
\eeqn
with the associated semigroup $S_{\Lambda_\eps}$, acting on the space of bounded Radon measures
$$
\XX:=M^1(\mathbb{R}_+) = \{ g \in (C_0(\R))'; \,\, \hbox{supp} \, g \subset \R_+ \},
$$
endowed with the weak $*$ topology $\sigma(M^1,C_0)$, where $C_0$
represents the space of continuous functions converging to $0$ at
infinity. From the duality of $S_{\Lambda^*_\eps}$, we deduce $S_{\Lambda_\eps}|_X=S_{L_\eps}$. The spectral analysis theory referred to \cite{KatoBook,Mbook*} 
indicates the structure of the spectrum denoting by $\Sigma(\Lambda_\eps)$ and the associated semigroup $S_{\Lambda_\eps}$.

\begin{theo}\label{th:MRe}
Assume \eqref{hyp:a1}-\eqref{hyp:a2}-\eqref{hyp:a3}-\eqref{hyp:a4} and
define  $\alpha := - a_0/2 < 0$. The operator $\Lambda_\eps$ is the
generator of a weakly $*$ continuous semigroup $S_{\Lambda_\eps}$ acting on $\XX$.
Moreover,  there exists a finite rank projector $\Pi_{\Lambda_\eps,\alpha}$ which commutes with $S_{\Lambda_\eps}$, an integer $j\ge0$ and some complex numbers
$$
\xi_1, ..., \xi_j \in \Delta_\alpha := \{ z \in \C, \,\, \Re e \, z > \alpha \},
$$
such that on $E_1 := \Pi_{\Lambda_\eps,\alpha}\XX$ the restricted operator satisfies
$$
\Sigma(\Lambda_{\eps|E_1} )\cap\Delta_\alpha=\{\xi_1,...,\xi_j\}
$$
(with the convention $\Sigma(\Lambda_{\eps|E_1} )\cap\Delta_\alpha = \emptyset$ when $j=0$) and for any $a > \alpha$ there exists a constant $C_a$ such that the remainder semigroup satisfies
$$
\|S_{\Lambda_\eps}(I-\Pi_{\Lambda_\eps,\alpha})\|_{\mathscr{B}(\XX)}\leq
C_a e^{a t}, \quad \forall \, t \ge 0.
$$
\end{theo}

In order to apply the Spectral Mapping Theorem of \cite{MS,Mbook*} and the
Weyl's Theorem of \cite{Voigt80,GMM*,MS,Mbook*}, we split the operator
$\Lambda_\eps$ as $\Lambda_\eps = \AA_\eps + \BB_\eps$ defined on $\XX$ by
\bear\label{def:Aeps}
   && \AA_\eps g:= \mu_\eps \MM_\eps[g], \quad \mu_\eps := \delta_0 -a'_\eps F_\eps, \\ \label{def:Beps}
   && \BB _\eps g:=-\partial_x g-a_\eps g.
\eear
As in the weak connectivity regime in \cite{MSWQ}, the properties
of the two auxiliary operators still hold in the strong one, which
implies Theorem \ref{th:MRe}.

\begin{lem}\label{lem:Lm1}
Assume that $a$ satisfies \eqref{hyp:a1}-\eqref{hyp:a2}-\eqref{hyp:a4}, then
the operators $\AA_\eps$ and $\BB _\eps$ satisfy the following properties.
\begin{enumerate}
  \item[(i)] $\AA_\eps\in\mathscr{B}(W^{-1,1}(\R_+),\XX)$.
  \item[(ii)] $S_{\BB _\eps}$ is $\alpha$-hypodissipative both in $X$ and $\XX$.
  \item[(iii)] The family of operators
  $S_{\BB _\eps}\ast\AA_\eps S_{\BB _\eps}$ satisfies
  $$
  \|(S_{\BB _\eps}\ast\AA_\eps S_{\BB _\eps})(t)\|_{\XX\rightarrow Y}\leq C_a e^{a t},\quad\forall a>\alpha.
  $$
for some positive constant $C_a$, where $Y:=BV(\R_+)\cap L^1_1(\R_+)$ with $BV(\R_+)$ representing the space of bounded variation measures while $L^1_1(\R_+)$ denoting the Lebesgue space weighted by $\langle x\rangle:=(1+|x|^2)^{1/2}$.
\end{enumerate}
\end{lem}

\proof (i) Under the assumption \eqref{hyp:a3} and \eqref{hyp:a4}, there exists some constant $K>0$ such that $|\eps\partial_\mu a|<K$, for any $\eps\ge0$, and we naturally have $\MM_\eps[\cdot]\in\BBB(W^{-1,1}(\R_+),\R)$, which implies that $\mu_\eps\in\BBB(\R,\XX)$. Then, we conclude as $\AA_\eps=\mu_\eps\MM_\eps[\cdot]$.

\smallskip
(ii) Here and below, we mark $A(x,\eps M_\eps)$ as $A(x)$ without obscurity. The explicit formula gives
$$
S_{\BB_\eps}(t)g(x)=e^{A(x-t)-A(x)}g(x-t)\mathbf{1}_{x-t\ge0}.
$$
From the assumption \eqref{hyp:a2}, there exists $x_1\in [0,\infty)$ such that $A(x)\ge{3\over 4}a_0(x-x_1)_+$, which implies
\beqn\label{eq:EoA1}
A(x-t)-A(x)\le C\,e^{3\beta t}, \quad 0\le t\le x,
\eeqn
where $C=e^{3a_0 x_1/4}$ and $\beta=-a_0/4<0$. Next, we conclude
$$
\|S_{\BB_\eps}(t)g\|_X\le C\,e^{3\beta t}\|g\|_X, \quad t\ge0,
$$
with  From (weakly *) density argument, we also have the same estimate in $\XX$. Then, $S_{\BB_\eps}$ is $\alpha-$ hypodissipative both in $X$ and $\XX$, as $\alpha>3\beta$. 

\smallskip
(iii) We denote
$$
N(t):=\MM_\eps[S_{\BB_\eps}(t)g]=(1-\kappa)^{-1}\int_0^\infty a_\eps e^{A(x-t)-A(x)}g(x-t)\mathbf{1}_{x-t\ge0}\mathrm{d}x.
$$
From the assumption \eqref{hyp:a3}, $N\in C^1_b(\R_+)$, we compute
\bean
N'(t) &=&  (1-\kappa)^{-1}\int_0^\infty\partial_x\Big(a_\eps e^{-A(x)}\Big)e^{A(x-t)}g(x-t)\mathbf{1}_{x-t\ge0}\mathrm{d}x\\
 &=& (1-\kappa)^{-1}\int_0^\infty(a'_\eps-a^2_\eps) e^{A(x-t)-A(x)}g(x-t)\mathbf{1}_{x-t\ge0}\mathrm{d}x.
\eean
The inequality \eqref{eq:EoA1} together with the assumption \eqref{hyp:a3} and \eqref{hyp:a4} permit us to have the following estimates
\bear
&|N(t)|\le(1-\kappa)^{-1}C a_1\ds\int_0^\infty e^{3\beta t}g(x)\mathrm{d}x\lesssim e^{3\beta t}\|g\|_\XX,&\label{ineq:EoN}\\
&|N'(t)|\le(1-\kappa)^{-1}C(K+a^2_1)\ds\int_0^\infty e^{3\beta t}g(x)\mathrm{d}x\lesssim e^{3\beta t}\|g\|_\XX.&\label{ineq:EoN'}
\eear
We continue to analyse the operator,
\bean
(S_{\BB _\eps}\ast\AA_\eps S_{\BB _\eps})(t)g(x) &=& \int_0^t(S_{\BB_\eps}(s)\mu_\eps)(x)N(t-s)\mathrm{d}s\\
 &=& \int_0^t e^{A(x-s)-A(x)}\mu_\eps(x-s)N(t-s)\mathbf{1}_{t-s\ge0}\\
 &=& e^{-A(x)}(\nu_\eps*\check{N}_t)(x),
\eean
where $\nu_\eps=e^{A}\mu_\eps$ and $\check{N}_t(\cdot)=N(t-\cdot)$. And its partial derivative with respect to $x$ is
\bean
\partial_x(S_{\BB _\eps}\ast\AA_\eps S_{\BB _\eps})(t)g &=& -a_\eps e^{-A(x)}(\nu_\eps*\check{N}_t)(x)-e^{-A(x)}(\nu_\eps*\check{N}'_t)(x)\\
 &&-e^{-A(x)}\nu_\eps(x-t)N(0)\mathbf{1}_{x-t\ge0}+e^{-A(x)}\nu_\eps(x)N(t).
\eean
Using the inequality \eqref{eq:EoA1},
\bean
\|e^{-A(x)}(\nu_\eps*\check{N}_t)(x)\|_{L^1_1} &\le& \|e^{\beta x}\langle x\rangle\|_{L^1}\|e^{-A(x)}(\nu_\eps*\check{N}_t)e^{-\beta x}\|_{L^\infty}\\
 &\lesssim& \Big\|\int_0^t e^{3\beta s}|\mu_\eps(x-s)N(t-s)|e^{-\beta x}\mathrm{d}s\Big\|_{L^\infty}\\
 &\lesssim& e^{2\beta t}\|(|\mu_\eps|e^{-\beta})*\check{(|N|e^{-2\beta})}_t\|_{L^\infty}\\
 &\lesssim& e^{\alpha t}\|\mu_\eps e^{-\beta}\|_\XX \|Ne^{-2\beta}\|_{L^\infty},
\eean
since $\|e^{\beta x}\langle x\rangle\|_{L^1}<\infty$. Then from the assumption \eqref{hyp:a4} as well as the estimates \eqref{ineq:EoF} and \eqref{ineq:EoN}, we get
\beqn\label{ineq:SBASB}
\|S_{\BB _\eps}\ast\AA_\eps S_{\BB _\eps})(t)g\|_{L^1_1} \lesssim e^{\alpha t}(1+K\|F_\eps e^{-\beta}\|_\XX) \|e^{\beta}\|_{L^\infty}\|g\|_\XX\lesssim e^{\alpha t}\|g\|_\XX.
\eeqn
Similarly, from \eqref{ineq:EoN'}, we have
\bean
\|e^{-A(x)}(\nu_\eps*\check{N}'_t)(x)\|_\XX &\lesssim& \int_0^t\int_0^\infty e^{3\beta s}|\mu_\eps(x-s)||N'(t-s)|\mathrm{d}s\\
 &\lesssim& e^{\alpha t}\|\mu_\eps\|_\XX\int_0^t e^{\beta s}\|g\|_\XX\mathrm{d}s\\
 &\lesssim& e^{\alpha t}\|g\|_\XX,
\eean
which implies
\bean
\|\partial_x(S_{\BB _\eps}\ast\AA_\eps S_{\BB _\eps})(t)g\|_\XX &\lesssim& a_1\|e^{-A}(\nu_\eps*\check{N}_t)\|_\XX+\|e^{-A}(\nu_\eps*\check{N}'_t)\|_\XX\\
&&+e^{3\beta t}\|e^{3\beta}\mu_\eps\|_\XX\|g\|_\XX+e^{3\beta t}\|e^{-A}\nu_\eps\|_\XX\|g\|_\XX\\
&\lesssim& e^{\alpha t}\|g\|_\XX.
\eean
We finally conclude $\mathrm{(iii)}$ from the above estimate together with \eqref{ineq:SBASB}.
\endproof


\subsection{Strong connectivity regime - exponential stability of the
  linearized equation}

Under the assumption \eqref{hyp:a4}, when the network connectivity parameter $\eps$ goes to infinity, the linearized time elapsed operator is simplified as
\begin{equation}\label{eq:TLC}
    \Lambda_\infty g=-\partial_x g-a(x,\infty)g+\delta_{x=0}\MM_\infty[g],
\end{equation}
where $\MM_\infty[g]=\ds\int_0^\infty a(x,\infty)g(x)\mathrm{d}x$.
Similarly to the vanishing limited case $\eps=0$ in  \cite{MSWQ}, we also obtain the following properties.
\begin{lem}\label{lem:Linfty}
In the limited case $\eps=\infty$, the operator $\Lambda_\infty$ and the associated semigroup $S_{\Lambda_\infty}$ satisfy
 \begin{enumerate}
  \item[(i)] $S_{\Lambda_\infty}$ is positive, i.e.
  $$
  S_{\Lambda_\infty}(t)g\in\XX_+,\quad\forall g\in\XX_+,\quad\forall t\ge0.
  $$
  \item[(ii)] $-\Lambda_\infty$ possesses the strong maximum principle, i.e. for any given $\mu\in\R$ and any nontrivial $g\in D(\Lambda_\infty)\cap\XX_+$, there holds
  $$
  (-\Lambda_\infty+\mu)g\ge0 \ implies \ g>0.
  $$
  \item[(iii)] $\Lambda_\infty$ satisfies the complex Kato's inequality, i.e.
  $$
  \Re e(\sgn g)\Lambda_\infty g\le\Lambda_\infty|g|, \quad\forall g\in D(\Lambda_\infty^2).
  $$
 \end{enumerate}
\end{lem}
Then we conclude the following exponential asymptotic stability.
\begin{theo}\label{th:lc}
There exist some constants $\alpha<0$ and $C>0$ such that $\Sigma(\Lambda_\infty) \cap \Delta_\alpha = \{ 0 \}$ and for any $g_0 \in X$, $\langle g_0 \rangle = 0$, there holds
\beqn\label{eq:lc}
\| S_{\Lambda_\infty}(t) g_0 \|_X \leq C e^{\alpha t} \, \| g_0 \|_X, \quad \forall \, t \ge 0.
\eeqn
\end{theo}
\proof Theorem \ref{th:SS} shows that there exists at least one nontrivial $F_\infty\ge0$ as the eigenvector to $0$ and the associated dual eigenvector is $\psi=1$. From Lemma \ref{lem:Linfty}-(ii)$\&$(iii), we deduce that the eigenvalue $0$ is simple and the associated eigenspace is Vect($F_\infty$). While Lemma \ref{lem:Linfty}-(i)$\&$(ii) imply that $0$ is the only eigenvalue with nonnegative real part. We then conclude from Theorem \ref{th:MRe}.
\endproof

We extend the exponential stability property in the limited case to the strong connectivity regime through a perturbation argument.

\begin{theo}\label{th:CWR}
There exist some constants $\eps_1 > 0$, $\alpha<0$ and $C>0$ such that for any $\eps \in [\eps_1,\infty]$ there hold $\Sigma(\Lambda_\eps) \cap \Delta_\alpha = \{ 0 \}$ and
\beqn\label{eq:CWR}
\| S_{\Lambda_\eps}(t) g_0 \|_X \leq C e^{\alpha t} \, \| g_0 \|_X, \quad \forall \, t \ge 0,
\eeqn
for any $g_0 \in X$, $\langle g_0 \rangle = 0$.
\end{theo}

The proof uses the stability theory for semigroups developed in Kato's book \cite{KatoBook} and revisited in \cite{MMcmp,T*,Mbook*}.
Now, we present several results needed in the proof of Theorem~\ref{th:CWR}.

\proof \noindent{\sl Step 1. Continuity of the operator.}
Directly from the definitions \eqref{def:Nepsg},  \eqref{def:Aeps} and \eqref{def:Beps} of $\MM_\eps$, $\AA_\eps$  and $\BB_\eps$, we have
$$
(\BB_\eps - \BB_\infty )g = (a(x,\infty)-a(x,\eps M_\eps ))g
$$
and
$$
(\AA_\eps - \AA_\infty )g = (\MM_\eps[g] - \MM_\infty[g]) \, \delta_0 -  \eps (\partial_\mu a)(x,\eps M_\eps)  \, F_\eps \, \MM_\eps[g].
$$
From the decay assumption \eqref{hyp:a4}, there
exists positive $\zeta_\eps\to 0$,  as $\eps\to+\infty$, such that
$|\eps\partial_\mu a(x,\eps M_\eps)|<\zeta_\eps$, for $\eps$ large enough. Together with the
smoothness assumption \eqref{hyp:a3}, we deduce that
\beqn\label{ineq:Leps-Linfty}
\| \BB_\eps - \BB_\infty \|_{\BBB(X)} + \| \AA_\eps - \AA_\infty \|_{\BBB(X)} \le C \, \zeta_\eps,
\eeqn
in the strong connectivity regime.

\smallskip
\noindent{\sl Step 2. Perturbation argument.} Similarly to \cite{MSWQ}, we present the sketch as the argument in the proof of \cite[Theorem 2.15]{T*} (see also \cite{KatoBook,MMcmp,Mbook*}). Define 
$$
K_\eps(z):=(R_{\BB_\eps} (z) \AA_\eps)^2 R_{\Lambda_\infty}(z)(\Lambda_\eps - \Lambda_\infty) \in \BBB(\XX,X).
$$
From Lemma \ref{lem:Lm1}-(i)$\&$(ii) as well as \eqref{ineq:Leps-Linfty}, there exist $\eps_1>0$ large enough and $C>0$, such that for any $z\in\Delta_\alpha\backslash B(0,\eta)$, for some $0<\eta<|\alpha|$ and any $\eps\in[\eps_1,\infty]$, the operator $K_\eps(z)$ satisfies
\beqn\label{ineq:Keps}
\| K_\eps (z) \|_{\BBB(X)} \le C\zeta_\eps\le C\zeta_{\eps_1}<1,
\eeqn
which permits us to well define $(1-K_\eps(z))^{-1}$ in $\BBB(X)$. From the Duhamel formula and the inverse Laplace transform, we have
$$
R_{\Lambda_\eps}=R_{\BB_\eps}-R_{\BB_\eps} \AA_\eps R_{\BB_\eps}+(R_{\BB_\eps} \AA_\eps)^2 R_{\Lambda_\eps}. 
$$
From the definition of $K_\eps$, we directly get
$$
 (I - K_\eps)R_{\Lambda_\eps}=R_{\BB_\eps}-R_{\BB_\eps} \AA_\eps R_{\BB_\eps}+ (R_{\BB_\eps} \AA_\eps )^2 R_{\Lambda_\infty}.
$$
From \eqref{ineq:Keps} and since for any $\eps \in [\eps_1,\infty]$, all the terms in the RHS of the above expression are clearly uniformly bounded in $\BBB(\XX,X)$ on $\Delta_\alpha \backslash B(0,\eta)$, we deduce that
$$
\Sigma(\Lambda_\eps) \cap \Delta_\alpha \subset B(0,\eta).
$$

\smallskip
Thanks to the unique continuity principle for holomorphic functions, we deduce that $R_{\Lambda_\eps} (b)|_X = R_{L_\eps} (b)$ for $b\in\R$ large enough, which implies that $R_{\Lambda_\eps} (z)|_X = R_{L_\eps} (z)$  for any $\Delta_\alpha \backslash B(0,\eta)$. By mean of Dunford integral (see \cite[Section III.6.4]{KatoBook} or \cite{GMM*,Mbook*}), we express the eigenprojector $\Pi_\eps$ as
\bean
\Pi_\eps &=& {i \over 2\pi} \int_{|z| = \eta} R_{\Lambda_\eps}(z)\,\mathrm{d}z
\\
&=& {i \over 2\pi} \int_{|z| = \eta} (I-K_\eps) R_{\Lambda_\eps}\,\mathrm{d}z+{i \over 2\pi} \int_{|z| = \eta}K_\eps R_{\Lambda_\eps}\,\mathrm{d}z
\\
&=& {i \over 2\pi} \int_{|z| = \eta}   (R_{\BB_\eps} \AA_\eps )^2 R_{\Lambda_\infty} \, dz
+  {i \over 2\pi} \int_{|z| = \eta}  K_\eps R_{\Lambda_\eps}  \, \mathrm{d}z,
\eean
where the contribution of holomorphic functions vanish. In a similar way, we have 
$$
\Pi_\infty={i \over 2\pi} \int_{|z| = \eta} R_{\Lambda_\infty}(z)\,\mathrm{d}z={i \over 2\pi} \int_{|z| = \eta}   (R_{\BB_\infty} \AA_\infty)^2 R_{\Lambda_0}\,\mathrm{d}z. 
$$
Next, we compute that
\bean
(R_{\BB_\eps} \AA_\eps )^2-(R_{\BB_\infty}\AA_\infty)^2
 &= & R_{\BB_\eps} \AA_\eps R_{\BB_\eps}  \{ (\AA_\eps-\AA_\infty) + (\BB_\infty-\BB_\eps) R_{\BB_\infty}\AA_\infty\}\\
 &&+ R_{\BB_\eps}  \{ (\AA_\eps- \AA_\infty)+\AA_\infty(\BB_\infty-\BB_\eps) R_{\BB_\infty}\} R_{\BB_\infty} \AA_\infty.
\eean
From the above identity together with the estimates \eqref{ineq:Leps-Linfty} and \eqref{ineq:Keps}, we deduce that
\bean
\| (\Pi_\eps - \Pi_\infty) g \|_{X} &=& \| (\Pi_\eps - \Pi_\infty) g \|_{\XX} \\
 &\le& {1 \over 2\pi}  \int_{|z| = \eta}   \| ((R_{\BB_\eps} \AA_\eps)^2  -  (R_{\BB_\infty} \AA_\infty)^2) R_{\Lambda_\infty} g \|_\XX\,\mathrm{d}z\\
 && +  {1 \over 2\pi}  \int_{|z| = \eta}   \| K_\eps R_{\Lambda_\eps} g \|_\XX \, \mathrm{d}z\\
 &\le& C \, \zeta_\eps \, \| g \|_\XX,
\eean
for any $g\in X$. Therefore, the eigenprojector $\Pi_\eps$ satisfies that
$$
\| \Pi_\eps - \Pi_\infty \|_{\BBB(X)} < 1,\quad\forall\eps\in[\eps_1,\infty].
$$

\smallskip
\noindent{\sl Step 3. Spectral gap.} From the classical result \cite[Section I.4.6]{KatoBook} (or more explicitly \cite[Lemma~2.18]{T*}), we deduce that there exists a unique simple eigenvalue $\xi_\eps \in \Delta_\alpha$ such that
$$
\Sigma(\Lambda_\eps)\cap\Delta_\alpha=\{\xi_\eps\},
$$
for any $\eps \in [\eps_1,+\infty]$. From the fact of mass conservation, we have $1 \in\XX'$ and $\Lambda_\eps^* 1 = 0$, which implies the conclusion since $\xi_\eps = 0$.
\endproof

\subsection{Strong connectivity regime - nonlinear exponential stability}

Now, we come back to the nonlinear problem \eqref{eq:ASM} in the case without delay as
$$
m(t)=p(t)=\int_0^\infty a(x,\eps m(t))f(x)\,\mathrm{d}x.
$$
In order to show that $m(t)$ is well defined, we recall the optimal transportation Monge-Kantorovich-Wasserstein distance on the probability measures set $\Ppp(\R_+)$ associated to the distance $d(x,y)= |x-y| \wedge 1$, denoting as $W_1$ and given by
$$
\forall \, f,g\in \Ppp(\R_+), \quad W_1(f,g) := \sup_{\varphi, \| \varphi \|_{W^{1,\infty}} \le 1} \int_0^\infty (f-g) \, \varphi.
$$
In addition, we define $\Phi :  L^1(\R_+) \times \R \to \R$ as
$$
\Phi [g,\mu] := \int_0^\infty a(x,\eps \mu) g (x) \, \mathrm{d}x - \mu.
$$

\begin{lem}\label{lem:varphig} Under the assumption \eqref{hyp:a2}-\eqref{hyp:a3}-\eqref{hyp:a4}, there exists $\eps_1 > 0$ large enough such that for any $\eps \in (\eps_1,+\infty)$, the equation $\Phi(g,\mu) = 0$ towards $\mu$ has a unique nonnegative solution $\mu = \varphi_\eps[g]$, where $\varphi_\eps : \Ppp(\R) \to \R$ is Lipschitz continuous for the weak topology of probability measures.
\end{lem}

\proof {\sl Step 1. Existence. }
For any $g \in \Ppp(\R)$, we obviously have $\Phi(g,0) > 0$ while for any $g \in \Ppp(\R)$ and $\mu>a_1$, we have
$$
 \Phi(g,\mu) \le a_1  -  \mu<0.
$$
Thanks to the intermediate value theorem, for any fixed $g \in \Ppp(\R_+)$ and any $\eps \ge 0$, there exists at least one solution $\mu \in (0,a_1]$ to the equation $ \Phi(g,\mu) = 0$ from the continuity property of $\Phi$.

\smallskip\noindent {\sl Step 2. Uniqueness and Lipschitz continuity. }
For any $f,g \in \Ppp(\R_+)$, from Step 1, we are able to consider $\mu,\nu \in \R_+$ such that
$$
\Phi(f,\mu) =  \Phi(g,\nu) = 0,
$$
which implies that
$$
\nu - \mu = \int_0^\infty a(x,\eps \nu) (g - f) +   \int_0^\infty (a(x,\eps \nu) -  a(x,\eps \mu) ) f.
$$
From the definition of $W_1$ and the assumption \eqref{hyp:a4}, we have
\bean
\Bigl| \int_0^\infty a(x,\eps \nu) (g - f) \Bigr| &\le& \| a (\cdot,\eps\nu) \|_{W^{1,\infty}} \, W_1(g,f),\\
\Bigl| \int_0^\infty \big(a(x,\eps \nu) -  a(x,\eps \mu)\big) f \Bigr| &\le& \| a(\cdot ,\eps \nu) -  a(\cdot,\eps \mu) \|_{L^{\infty}}\le\zeta_\eps |\mu-\nu|,
\eean
where $\zeta_\eps\to 0$ as $\eps\to\infty$. We then take $\eps_1$ large enough, such that
$$
1-\zeta_\eps\in(0,1),\quad\forall\eps\in[\eps_1,\infty],
$$
which permitting us to obtain
\beqn\label{eq:mu-nu}
|\mu - \nu| \, (1-\zeta_\eps) \le \| a (\cdot,\eps\nu) \|_{W^{1,\infty}} \, W_1(g,f).
\eeqn
On the one hand, when $f=g$, we immediately deduce that $\mu=\nu$, the uniqueness of the definition of the mapping $\varphi_\eps[g] := \mu$. On the other hand, we get the Lipschitz continuity of the function $\varphi_\eps$ directly from the inequation \eqref{eq:mu-nu}.
\endproof

We present the proof of our main result Theorem~\ref{th:MR} in the case without delay.

\proof[Proof of Theorem~\ref{th:MR} in the case without delay] We split the proof into three steps.

\noindent
{\sl Step 1. New formulation.} Benefiting from Lemma \ref{lem:varphig}, in the strong connectivity regime $\eps \in [\eps_1,\infty)$, where $\eps_1$ is the same as that in Lemma \ref{lem:varphig}, we introduce a new formulation of the solution $f \in C([0,\infty); X)$ to the evolution equation \eqref{eq:ASM} and the solution $F_\eps$ to the stationary problem \eqref{eq:StSt} satisfying
\bean
\partial_t f + \partial_x f +a(\eps \varphi[f]) f=0, &\quad& f(t,0) =  \varphi[f(t,\cdot)],\\
\partial_x F +a(\eps M) F=0, &\quad&  F (0) =  M = \varphi[F],
\eean
for a given unit mass initial datum $0 \le f_0 \in X$, where here and below the $\eps$ and $x$ dependency is often removed without any confusion.

\smallskip
Next, we consider the variation function $g:= f - F$ satisfying
\bear
\partial_t g=-\partial_x g - a(\eps M) g - \eps a'(\eps M) F \, \MM[g]-Q[g],
\eear
where
$$
Q[g] := a(\eps \varphi[f]) f - a(\eps \varphi[F]) F  - a(\eps
\varphi[F]) g - \eps a'(\eps  \varphi[F]) F \MM[g],
$$
with $\MM = \MM_\eps$ defined in \eqref{def:Nepsg}, complemented with the boundary condition given by
\bean
g(t,0) &=&  \varphi[f(t,\cdot)] - \varphi[F]\\
 &=& \int_0^\infty a(\eps \varphi[f]) f  - \int_0^\infty a(\eps \varphi[F]) F\\
&=&  \MM[g] +  \QQ[g],
\eean
where $\QQ[g] := \langle Q[g] \rangle$. Regarding the boundary condition as a source term again, we deduce that the variation function $g$ satisfies the equation
\beqn\label{eq:dtg=Lambda+Z}
\partial_t g = \Lambda_\eps g + Z[g],
\eeqn
with the nonlinear term $Z[g] :=  - Q[g] + \delta_0 \QQ[g]$.

\smallskip
\noindent{\sl Step 2. The nonlinear term.} With the fact that $f$ is mass
conserved, $\|F\|_X=1$ and the assumption \eqref{hyp:a4}, we estimate that
\bean
\|Q[g]\|_X &=& \|a(\eps\varphi[f])f-a(\eps\varphi[F])f-\eps
a'(\eps\varphi[F])F\,\MM[g]\|_X\\
&\le& \eps\|a'\|_{L^\infty_x}\|f\|_X\big|\varphi[f]-\varphi[F]\big|+\eps\|a'\|_{L^\infty_x}\|F\|_X\MM[g]\\
&\lesssim& \zeta_\eps\big(\MM[g]+\|Q[g]\|_X\big)+\zeta_\eps\MM[g],
\eean
where $\zeta_\eps\to0$ as $\eps\to+\infty$. Considering that
$$
\MM[g]\le a_1(1-\kappa)^{-1}\|g\|X\lesssim\|g\|_X,
$$
from the above inequality, we deduce that
$$
\|Q[g]\|\lesssim\zeta_\eps\MM[g]\lesssim\zeta_\eps\|g\|_X,
$$
with $\eps$ large enough. We then obtain
$$
\|Z[g]\|_X\le 2\|Q[g]\|_X\lesssim\zeta_\eps\|g\|_X
$$

\smallskip
\noindent{\sl Step 3. Decay estimate.} Thanks to the Duhamel formula, the solution $g$ to the evolution equation \eqref{eq:dtg=Lambda+Z}
satisfies
$$
g(t) = S_{\Lambda_\eps}(t) (g_0) + \int_0^t S_{\Lambda_\eps}(t-s)  Z[g(s)] \, \mathrm{d}s.
$$
Benefiting from Theorem~\ref{th:CWR} and the second step, we deduce
\bean
\| g(t) \|_X
 &\le&  C \, e^{\alpha t} \, \| g_0 \|_X + \int_0^t C \, e^{\alpha (t-s) } \, \| Z[g(s)] \|_X  \, \mathrm{d}s
\\
 &\lesssim&  e^{\alpha t} \, \| g_0 \|_X + \zeta_\eps\int_0^t  e^{\alpha (t-s) } \, \| g(s) \|_X  \, \mathrm{d}s,
 \eean
for any $t \ge 0$ and for some constant $\alpha < 0$, independent of
$\eps \in [\eps_1,+\infty)$. Thanks to the Gronwall's lemma (for linear integral inequality), we have
\bean
\|g(t)\|_X &\lesssim& e^{\alpha t}\|g_0\|_X+\zeta_\eps\|g_0\|_X\int_0^t
e^{\alpha t}\exp\{\int_s^t e^{\alpha(t-r)}\mathrm{d}r\}\mathrm{d}s\\
&\lesssim& e^{\alpha t}\|g_0\|_X+\zeta_\eps te^{\alpha t}\|g_0\|_X\\
&\lesssim& e^{\alpha' t}\|g_0\|_X,
\eean
for some constant $\alpha<\alpha'<0$.
\endproof

\section{Case with delay}
\label{sec:WithDelay}

In this section, we conclude our main result Theorem~\ref{th:MR}, in the case with delay by following the same strategy as section \ref{sec:WithoutDelay} but with appropriate adaptation towards the boundary term. Recalling from Theorem~\ref{th:SS}, we already know that their exists a uinque stationary state $(F_\eps, M_\eps)$ in the strong connectivity regime, therefore, we start from the linearization of the evolution equation \eqref{eq:ASM}.

\subsection{Linearized equation and structure of the spectrum} We still consider the variation functions $(g,n,q)=(f,m,p)-(F_\eps,M_\eps,M_\eps)$ around the steady state. However, because of the failure to express $n(t)$ explicitly, we introduce the following intermediate evolution equation on a function $v=v(t,y)$
\begin{equation}\label{eq:Eov}
    \partial_t v+\partial_y v=0, \ \
    v(t,0)=q(t), \ \ v(0,y)=0,
\end{equation}
where $y\geq0$ represent the local time for the network activity. Clearly, the above equation can be solved by the characteristics method as
$$
v(t,y)=q(t-y)\mathbf{1}_{0\leq y\leq t},
$$
which simplifies the expression of the variation $n(t)$ of network activity in the equation \eqref{eq:ASMlin3}, given by
$$
n(t)= \DD[v(t)], \quad \DD[v] := \int_0^\infty v(y)b(\mathrm{d}y),
$$
while the variation $q(t)$ of discharging neurons in the equation \eqref{eq:ASMlin2} is simplified as
$$
q(t)= \OO_\eps[g(t),v(t)]:= \NN_\eps[g(t)] + \kappa_\eps  \, \DD[v(t)],
$$
where
$$
\NN_\eps[g]  := \int_0^\infty a_\eps(M_\eps) g \, \mathrm{d}x, \quad \kappa_\eps := \int_0^\infty a'_\eps (M_\eps) F_\eps \, \mathrm{d}x.
$$
Therefore, we may rewrite the linearized system \eqref{eq:ASMlin1}-\eqref{eq:ASMlin2}-\eqref{eq:ASMlin3} together with the evolution equation \eqref{eq:Eov} as an autonomous system
\begin{equation}\label{eq:AS}
\partial_t
\begin{pmatrix}
g\\v
\end{pmatrix}
=\LLL_\eps
\begin{pmatrix}
g\\v
\end{pmatrix}
:=
\begin{pmatrix}
-\partial_x g-a_\eps g-a'_\eps F_\eps \DD [v]\\
-\partial_y v 
\end{pmatrix},
\end{equation}
with the associated semigroup $S_{\LLL_\eps}(t)$ acting on
$$
X:=X_1\times X_2:=L^1(\R_+)\times L^1(\R_+,\mu),
$$
where the measure $\mu(x)=e^{-\delta x}$ with the same $\delta>0$ in the condition \eqref{hyp:del}.
The domain of the operator $\LLL_\eps$ is given by
$$
D(\LLL_\eps):=\{(g,v)\in W^{1,1}(\R_+)\times W^{1,1}(\R_+,\mu); \ g(0)=v(0)=\OO_\eps[g,v]\}.
$$
We also consider the boundary condition as a source term and rewrite the autonomous systemwe as
$$
\partial_t (g,v)=\varLambda_\eps(g,v),
$$
where the generator $\varLambda_\eps = (\varLambda^1_\eps,\varLambda^2_\eps)$ is given by 
\bean
\varLambda^1_\eps(g,v) &:=& -\partial_x g-a_\eps g-a'_\eps
F_\eps \DD [v]+\delta_{x=0}\OO_\eps[g,v],\\
\varLambda^2_\eps(g,v) &:=& -\partial_y v +\delta_{y=0}\OO_\eps[g,v].
\eean
and the associated semigroup $S_{\varLambda_\eps}(t)$ acting on
$$
\XX:=\XX_1\times \XX_2:=M^1(\mathbb{R}_+)\times M^1(\mathbb{R}_+,\mu).
$$
Similarly, we have $S_{\varLambda_\eps}|_X=S_{\LLL_\eps}$.  Next, we are going to show that $S_{\varLambda_\eps}$ also possesses a suitable decomposition of a finite dimensional principal part as well as an exponential decaying remainder.

\begin{theo}\label{th:MRe1} Under the assumptions \eqref{hyp:a1}-\eqref{hyp:a2}-\eqref{hyp:a3}-\eqref{hyp:a4} as well as the condition \eqref{hyp:del} and taking delay into account, the conclusions of Theorem~\ref{th:MRe} still holds true with $\alpha := \max\{-a_0/2,-\delta\}<0$.
\end{theo}

The result also comes from the spectral analysis approach. In order to apply the Spectral Mapping theorem and the Weyl's Theorem established in \cite{MS,Mbook*}, we split the operator appropriately as $\varLambda_\eps =\AA_\eps+\BB _\eps$ with
$$
\BB_\eps (g,v) = \begin{pmatrix}
\BB^1_\eps(g,v)\\ \BB^2_\eps(g,v)
\end{pmatrix}
=
\begin{pmatrix}
-\partial_x g-a_\eps g\\
-\partial_y v
\end{pmatrix}
$$
and
$$
\AA_\eps (g,v) =  \begin{pmatrix}
\AA^1_\eps(g,v)\\ \AA^2_\eps(g,v)
\end{pmatrix}
=
\begin{pmatrix}
-a'_\eps F_\eps\DD[v]+\delta_{x=0}\OO _\eps[g,v]\\
\delta_{y=0}\OO _\eps[g,v]
\end{pmatrix},
$$
which hold the following parallel properties as those in Lemma \ref{lem:Lm1}.
\begin{lem}\label{lem:Lm2}
\begin{enumerate}
  \item[(i)] $\AA_\eps\in\BBB(W^{-1,1}(\R_+)\times W^{-1,1}(\R_+,\mu),\XX)$.
  \item[(ii)] $S_{\BB_\eps}(t)$ is $\alpha$-hypodissipative in both $X$ and $\XX$;
  \item[(iii)] the family of operators
  $S_{\BB_\eps}\ast\AA_\eps S_{\BB_\eps}$ satisfies
  $$
  \|(S_{\BB_\eps}\ast\AA_\eps S_{\BB_\eps})(t)\|_{\BBB(\XX,Y)}\leq C_a e^{\alpha
  t},\quad\forall a>\alpha,\quad\forall t\ge0,
  $$
  for some constant $C_a>0$ and with $Y:=Y_1\times Y_2$, where $Y_1=BV(\R_+)\cap L_1^1(R_+)$ and $Y_2=BV(\R_+,\mu)\cap L_1^1(R_+,\mu)$.
\end{enumerate}
\end{lem}
We skip the proof and refer to the proof of Lemma 3.2 in \cite{MSWQ} for more details.

\subsection{Strong connectivity regime - exponential stability of the linearized equation}

In the limited case, i.e. as the network connectivity parameter $\eps$ passes to the infinity, $\kappa_\eps$ vanishes from the assumption \eqref{hyp:a4}, which simplifies the linearized operator as
\begin{equation}\label{eq:lc1}
    \varLambda_\infty
    \begin{pmatrix}
    g \\ v
    \end{pmatrix}
    =
    \begin{pmatrix}
    -\partial_x g-a(x,\infty)g+\delta_{x=0}\OO _\infty[g,v]\\
    -\partial_y v  +\delta_{y=0}\OO _\infty[g,v]
    \end{pmatrix},
\end{equation}
where $\OO _\infty[g,v]=\NN_\infty[g]=\ds\int_0^\infty
a(x,\infty)g(x)\mathrm{d}x$. We have already proven the exponential stability of the first component $\varLambda_\infty^1$,  then the Duhamel formula
$$
v(t)=S_{\BB^2_\infty}(t)v_0+\int_0^t S_{\BB^2_\infty}(t-s)\AA^2_\infty\Big(g(s),v(s)\Big)\mathrm{d}s
$$
implies the similar exponential asymptotic estimate for the second component $\varLambda_\infty^2$. Together with Theorem \ref{th:lc}, we have
\begin{theo}\label{th:lc1}
There exist some constants $\alpha<0$ and $C>0$ such that
$\Sigma(\varLambda_\infty) \cap \Delta_\alpha = \{ 0 \}$ and for any $(g_0, \, v_0) \in X$, $\langle g_0 \rangle = 0$, there holds
\beqn\label{eq:lc1}
\| S_{\varLambda_\infty}(t) (g_0, \,v_0) \|_\XX \leq C e^{\alpha t} \, \| (g_0, \, v_0) \|_\XX, \quad \forall \, t \ge 0.
\eeqn
\end{theo}
Then we extend the geometry structure of the
spectrum of the linearized time elapsed equation in the limited case to the strong
connectivity regime taking delay into account.
\begin{theo}\label{th:CWR1}
There exist some constants $\eps_1 > 0$, $C \ge 1$ and $\alpha < 0$ such that for any $\eps\in[\eps_1,+\infty]$
there holds $\Sigma(\varLambda_\eps)\cap\Delta_\alpha=\{0\}$ and
\beqn\label{eq:Cwd1}
\|S_{\varLambda_\eps}(t)(g_0,v_0)\|_\XX\leq C e^{\alpha t}\|(g_0,v_0)\|_\XX,
\eeqn
for any $(g_0,v_0) \in\XX$ such that $\langle g_0 \rangle = 0$.
\end{theo}

\proof We proceed in two steps.

\noindent{\sl Step 1. Continuity of the operator $\varLambda_\eps$.} For all $(g,v)\in\XX$, we have
\begin{subequations}
 \begin{align}
 &\varLambda_\eps
 \begin{pmatrix}
 g \\ v
 \end{pmatrix}
 =
 \begin{pmatrix}
 -\partial_x g-a_\eps g-a'_\eps
 F_\eps\DD_\eps[v]+\delta_{x=0}\OO_\eps[g,v]\\
 -\partial_y v +\delta_{y=0}\OO_\eps[g,v]
 \end{pmatrix},\label{eq:Le1}\\
 &\varLambda_\infty
 \begin{pmatrix}
 g \\ v
 \end{pmatrix}
 =
 \begin{pmatrix}
 -\partial_x g-a(x,\infty) g +\delta_{x=0}\OO_\infty[g,v]\\
 -\partial_y v +\delta_{y=0}\OO_\infty[g,v]
 \end{pmatrix}. \label{eq:L01}
 \end{align}
\end{subequations}
Compute the difference between \eqref{eq:Le1} and \eqref{eq:L01}, we have
\begin{equation*}
(\varLambda_\eps-\varLambda_\infty)
\begin{pmatrix}
g \\ v
\end{pmatrix}
=
\begin{pmatrix}
(a(x,\infty)-a_\eps)g-a'_\eps
F_\eps\DD_\eps[v]+\delta_{x=0}\big(\OO_\eps[g,v]-\OO_\infty[g,v]\big) \\
\delta_{y=0}\big(\OO_\eps[g,v]-\OO_\infty[g,v]\big)
\end{pmatrix}.
\end{equation*}
From the assumption \eqref{hyp:a2} and \eqref{hyp:a4}, we deduce that
\begin{eqnarray*}
\|(\varLambda_\eps-\varLambda_0)(g,v)\|_\XX &=& \|(a(\cdot,\infty)-a_\eps)g\|_{\XX_1}+\|a'_\eps
  F_\eps\DD_\eps[v]\|_{\XX_1}+2|\OO_\eps[g,v]-\OO _\infty[g,v]| \\
 &\leq& 3\|(a_\infty-a_\eps)g\|_{\XX_1}+2\|a'_\eps F_\eps\DD_\eps[v]\|_{\XX_2}\\
 &\leq& 3\zeta_\eps\|g\|_{\XX_1}+2a_1\zeta_\eps(1-\zeta_\eps)\|F_\eps\|_{\XX_1} \|v\|_{\XX_2}\\
 &\lesssim& \zeta_\eps\|(g,v)\|_\XX,
\end{eqnarray*}
where $\zeta_\eps\to0$ as $\eps\to\infty$, which implies the continuity of the operator$\varLambda_\eps$ to $\eps$ in the strong connectivity regime.

\noindent{\sl Step 2. Extension to the strong connectivity regime.}
Similar to the proof of Theorem~\ref{th:CWR}, we deduce that their exists $\eps_1>0$ large enough, such that for any $\eps\in[\eps_1,\infty]$, the eigenprojector $\Pi_\eps$ satisfies
$$
\|\Pi_\eps-\Pi_\infty\|_{\BBB(X)}<1,
$$
which permits us to conclude that there exists $\xi_\eps$ satisfying $|\xi_\eps|\leq O(\zeta_\eps)$, such that
$$
\Sigma(\varLambda_\eps)\cap\Delta_\alpha=\{\xi_\eps\}
$$
and $\xi_\eps$ is algebraically simple (see in \cite[Section I.4.6]{KatoBook} and \cite{T*}). Since the dual operator $\varLambda_\eps^*$ is given by
$$
\varLambda^*_\eps
\begin{pmatrix}
\varphi \\ \psi
\end{pmatrix}
=
\begin{pmatrix}
\partial_x \varphi -a_\eps \varphi + a_\eps (\varphi(0) + \psi(0))
\\
\partial_y  \psi + \kappa_\eps b \, \psi(0) + \kappa_\eps b \, \varphi(0) - b\ds\int a_\eps' F_\eps \, \varphi \, \mathrm{d}x
\end{pmatrix},
$$
for any $(\varphi,\psi)\in\XX'$, observe that $\varLambda_\eps^*(1,0) = 0$. Thus, $0\in\Sigma(\varLambda_\eps^\ast)$, which implies $\xi_\eps=0$. Together with the fact that $\langle g_0 \rangle = \langle (g_0,v_0),(1,0) \rangle_{\XX,\XX'} =0$, the exponential asymptotic estimate \eqref{eq:Cwd1} holds.
\endproof

\subsection{Strong connectivity regime - nonlinear exponential stability}

Finally, we focus on the nonlinear problem taking delay into account and present the proof of the first part of our main result in the case with delay, neglecting the connectivity parameter $\eps$ most of time without any misleading.

\proof[Proof of Theorem~\ref{th:MR} in case with delay]
Inspired by the intermediate evolution equation \eqref{eq:Eov}, we rewrite the nonlinear problem as a system of $(f,u)$, given by
\bean
\partial_t f &=& - \partial_x f - a_\eps(\DD[u]) f + \delta_0 \PPP[f,\DD[u]]
 \\
 \partial_t u&=& - \partial_y u  + \delta_0 \PPP[f,\DD[u]],
\eean
where
$$
\PPP[f,m] = \int a(m) f, \quad \DD[u] = \int b u.
$$
Denoting $U:= M {\bf 1}_{y \ge 0}$, the steady state $(F,U)$ satisfies
\bean
0 &=& - \partial_x F- a_\eps(M) F + \delta_0 M\\
0 &=& - \partial_y U  + \delta_0 M, \quad M = \DD[U] = \PPP[F,\DD[U]].
\eean
Introducing the variation functions $g:= f - F$ and $v = u - U$ again, we obtain the system of $(g,v)$ as
\bean
\partial_t g
&=& - \partial_x g - a_\eps(\DD[u]) f + a_\eps(M) F + \delta_0 (\PPP[f,\DD[u]] - \PPP[F,\DD[U]])
\\
&=& - \partial_x g - a_\eps(M) g -a'_\eps F \DD[v] - Q[g,v] + \delta_0 \OO[g,v] + \delta_0 \QQ[g,v]
\\
&=& \varLambda_\eps^1(g,v) + \ZZ^1[g,v],\\
\partial_t v
&=& - \partial_y v  + \delta_0 (\PPP[f,\DD[u]] - \PPP[F,\DD[U]])
\\
&=& - \partial_y v +  \delta_0 \OO[g,v] + \delta_0 \QQ[g,v]
\\
&=& \varLambda_\eps^2(g,v) + \ZZ^2[g,v],
\eean
where
\bean
Q[g,v] &:=& a_\eps(\DD[u])f-a_\eps(M) F-a_\eps(M) g-a'_\eps F \DD[v]\\
\QQ[g,v] &:=& \langle Q[g,v] \rangle,
\eean
and the remainders are given by
\bean
\ZZ^1[g,v] &:=& - Q[g,v] + \delta_0 \QQ[g,v],\\
\ZZ^2[g,v] &:=& \delta_0 \QQ[g,v].
\eean
From the mass conservation of $f$ and the assumption \eqref{hyp:a4}, we deduce that
\bean
\|Q[g,v]\|_{X_1} &=& a_\eps(\DD[u])f-a_\eps(M)f-a'_\eps F\DD[v]\\
&\le& \zeta_\eps\|f\|_{X_1}\Big|\DD[u]-\DD[U]\Big|-\zeta_\eps\|\DD[v]\\
&\lesssim& \zeta_\eps\|v\|_{X_2},
\eean
where $\zeta_\eps$ is the same as in Theorem \ref{th:CWR1}. Denoting $\ZZ[g,v]:=(\ZZ^1[g,v],\ZZ^2[g,v])$, we clearly have the estimate
$$
\|\ZZ[g,v]\|_X\lesssim\zeta_\eps\|(g,v)\|_X
$$
The associated Duhamel formula writes
$$
(g(t),v(t)) = S_{\varLambda_\eps}(t) (g_0,v_0) + \int_0^t S_{\varLambda_\eps}(t-s) \ZZ[g(s),v(s)] \, \mathrm{d}s.
$$
Using the above estimate for the nonlinear term and the Gronwall's
lemma, we conclude as in the proof of Theorem \ref{th:MR}.
\endproof

\section{Step function firing rate}\label{sec:StepFctFR}
In this section, we focus on the nonlinear time elapsed model in
\cite{PK1,PK2} with a particular step function firing rate given by
$$
a(x,\mu)=\mathbf{1}_{x>\sigma(\eps\mu)}
$$
satisfying \eqref{hyp:s1}-\eqref{hyp:s2}-\eqref{hyp:s3}-\eqref{hyp:s4}. We consider the dynamic of the neuron network \eqref{eq:ASM} completed with an initial probability density $f_0$ satisfying
\begin{equation}\label{eq:ID}
0\leq f_0\leq 1, \quad \int_0^\infty f_0(x)\mathrm{d}x=1.
\end{equation}
Obviously, the solution $f$ of the time elapsed equation \eqref{eq:ASM} corresponding to the firing rate \eqref{hyp:s1} is still mass conserved,  and we naturally renormalize that mass. The previous work \cite{PK1} shows that the model \eqref{eq:ASM} with the step function firing rate \eqref{hyp:s1} admits a steady state as well as a unique solution, that is to say the second part of Theorem \ref{th:EaU} and Theorem \ref{th:SS}. By applying the adapted above spectral analysis method, we conclude the results in Theorem \ref{th:MR} for a particular step function firing rate, which accurate the stability results in \cite{PK1} in the case with delay. Failing to construct the linearized equations \eqref{eq:ASMlin1} and \eqref{eq:ASMlin2}, we replace them with another more concise linear equation for the variation functions $(g,n,q) = (f,m,p) - (F_\eps,M_\eps,M_\eps)$, which writes
\bean
&& \partial_t g+\partial_x g+g\mathbf{1}_{x>\sigma_\eps}=0,\\
&& g(t,0)=q(t):=\int_0^\infty g\mathbf{1}_{x>\sigma_\eps}\mathrm{d}x, \quad g(0,x)= g_0(x),
\eean
where here and below we note $\sigma_\eps:=\sigma(\eps M_\eps)$ for
simplicity. We introduce the intermediate evolution equation \eqref{eq:Eov} again
to write the linear equation \eqref{eq:ASMFRlin1}-\eqref{eq:ASMFRlin2}-
\eqref{eq:ASMlin3} as a time autonomous system
\begin{equation}\label{eq:ASFR}
\partial_t(g,v)=\LLL_\eps(g,v),
\end{equation}
where the operator $\LLL_\eps = (\LLL^1_\eps,\LLL^2_\eps)$ is defined by 
\bean
\LLL^1_\eps(g,v) &:=& -\partial_x g-g\mathbf{1}_{x>\sigma_\eps}+\delta_{x=0}\RR_\eps[g,v],
\\
\LLL^2_\eps(g,v) &:=& -\partial_y v +\delta_{y=0}\RR_\eps[g,v], 
\eean
with the boundary term
$$
\RR_\eps[g(t),v(t)]:= \PP[g,M_\eps] = \int_0^\infty g \mathbf{1}_{x>\sigma_\eps} \,\mathrm{d}x,
$$
in the space
$$
X=X_1\times X_2:=L^p_0(\R_+)\times  L^p(\R_+,\mu)
$$
with $L^p_0(\R_+)=\{h\in L^p(\R_+);\ \langle h\rangle=0\}$, $1\leq p\leq\infty$, and $\mu(x)=e^{-\delta x}$, $\delta>0$ is the same as in the condition \eqref{hyp:del}. We extend the exponential stability from the single equation of $g$ to the above autonomous system.

\begin{theo}\label{th:MReFR} Assume \eqref{hyp:s1}-\eqref{hyp:s2}-\eqref{hyp:s3} and \eqref{hyp:del} (with \eqref{hyp:s4}). There exist some constants $\eps_0 > 0$ $(\eps_1>0)$, $C \ge 1$ and
$\alpha < 0$ such that for any $\eps\in[0,\eps_0]$ $\eps\in[\eps_1,\infty]$
there holds $\Sigma(\LLL_\eps)\cap\Delta_\alpha=\{0\}$ and 
\beqn\label{CwdFR}
\|S_{\LLL_\eps}(t)(g_0,v_0)\|_X\leq C e^{\alpha t}\|(g_0,v_0)\|_X,
\eeqn
for any $(g_0,v_0) \in X$, s.t. $\langle (g_0,v_0)\, , (1,0)
\rangle_{X,X'} = 0$. 
\end{theo}
The extension follows from the Spectral Mapping theorem and the Weyl's Theorem
by introducing a convenient splitting of the operator $\LLL_\eps$ as
$\LLL_\eps =\AA_\eps+\BB _\eps$ with
$$
\BB_\eps (g,v) = \begin{pmatrix}
\BB^1_\eps(g,v)\\ \BB^2_\eps(g,v)
\end{pmatrix}
=
\begin{pmatrix}
-\partial_x g- g\mathbf{1}_{x>\sigma_\eps}\\
-\partial_y v 
\end{pmatrix}
$$
and
$$
\AA_\eps (g,v) =  \begin{pmatrix}
\AA^1_\eps(g,v)\\ \AA^2_\eps(g,v)
\end{pmatrix}
=
\begin{pmatrix}
\delta_{x=0}\RR _\eps[g,v]\\
\delta_{y=0}\RR _\eps[g,v]
\end{pmatrix}.
$$
Since the step function firing rate is no longer continuous, we have to consider the resolvent of the operator $\BB_\eps$.

\begin{lem}\label{lem:LmFR}Assume \eqref{hyp:s1}-\eqref{hyp:s2}-\eqref{hyp:s3} and \eqref{hyp:del}
  (with \eqref{hyp:s4}). Then the two operators satisfy
\begin{enumerate}
  \item[(i)] $\AA_\eps\in\mathscr{B}(X,Y)$, where
    $Y=\C\delta_0\times\C\delta_0\subset
    X$ with compact embedding;
  \item[(ii)] $S_{\BB_\eps}(t)$ is $\alpha$-hypodissipative in $X$;
  \item[(iii)] $(R_{\BB_\eps}(z)\AA_\eps)^2(z)\in\BBB(X)$, with bound in $\OO(\langle z\rangle^{-1})$, $\forall z\in\Delta_{-1}$.
  \end{enumerate}
\end{lem}

\proof In order to simplify the notation, we note
$$
\rho(x):=\int_0^x\mathbf{1}_{y>\sigma_\eps}\mathrm{d}y=(x-\sigma_\eps)_+.
$$
\smallskip
(i) It is an immediate consequence of the fact that $\DD \in \BBB(X_2;\R)$ (because of \eqref{hyp:del}) and $\NN_\eps \in \BBB(X_1;\R)$.\\
\smallskip
(ii) We write $S_{\BB^1_\eps}$ and $S_{\BB^2_\eps}$ respectively with the explicit formula
\bean
  S_{\BB^1_\eps}(t)g(x) &=& e^{\rho(x-t)-\rho(x)}g(x-t)\mathbf{1}_{x-t\geq0},\\
  S_{\BB^2_\eps}(t)v(y) &=& v(y-t)\mathbf{1}_{y-t\geq0}.
\eean
We estimate that 
  \bean
  \| S_{\BB^1_\eps}(t)g\|_{X_1} &=& \|e^{\rho(x)-\rho(x+t)}g(x)\|_{X_1}=\|e^{(x-\sigma_\eps)_+-(x+t-\sigma_\eps)_+}g(x)\|_{X_1}\\
  &\le& Ce^{-t}\|g(x)\|_{X_1},\\
  \| S_{\BB^2_\eps}(t) v \|_{X_2} &=& \|e^{-\delta(y+t)}v(y)\|_{L^p}=e^{-\delta t}\| v \|_{X_2}, 
  \eean
  for any $g\in X_1$ and $v\in X_2$ and any $t \ge 0$, which implies
  $$
  \|S_{\BB_\eps}(t)\|_{X\to X}\le C e^{\alpha t}, \quad t\ge 0,
  $$
  by choosing $C:=\max\{2e^{\sigma_\eps},1\}$ and $\alpha:=\max\{-1,-\delta\}$.\\  
\smallskip
(iii) We have 
$$
S_{\BB^1_\eps}(t)\AA_\eps[g,v](x)=\delta_{x=t}\RR_\eps[g,v]e^{\rho(x-t)-\rho(x)},
$$
and we denote
\bean
k_t(x) &:=&\AA_\eps S_{\BB^1_\eps}(t)\AA_\eps[g,v](x)=\delta_{x=0}\int_0^\infty\delta_{x=t}e^{\rho(x-t)-\rho(x)}\mathbf{1}_{x>\sigma_\eps}\mathrm{d}x\\
 &=& \delta_{x=0}e^{-\rho(t)}\mathbf{1}_{t>\sigma_\eps}\RR_\eps[g,v].
\eean
Finally, we obtain
\bean
(S_{\BB^1_\eps}\AA_\eps)^{(\ast2)}(t)[g,v](x) &=& \int_0^t k_{t-s}(x-s)e^{\rho(x-s)-\rho(x)}\mathbf{1}_{x-s\ge0}\mathrm{d}s\\
 &=& e^{-\rho(t-x)+\rho(0)-\rho(x)}\mathbf{1}_{t\ge x}\mathbf{1}_{t-x>\sigma_\eps}\RR_\eps[g,v]\\
 &=& e^{-\rho(t-x)-\rho(x)}\mathbf{1}_{t-x>\sigma_\eps}\RR_\eps[g,v].
\eean
Denoting $\psi_t(x):=e^{-\rho(t-x)-\rho(x)}\mathbf{1}_{t-x>\sigma_\eps}$, we compute its Laplace transform
\bean
\hat{\psi}(z) &=& \int_0^\infty e^{-\rho(t-x)-\rho(x)}\mathbf{1}_{t-x>\sigma_\eps}e^{-zt}\mathrm{d}t\\
 &=& e^{-\rho(x)-zx}\int_{\sigma_\eps}^\infty e^{(t-\sigma_\eps)_+-zt}\mathrm{d}t\\
 &=&  e^{-\rho(x)-z(x+\sigma_\eps)}\int_0^\infty e^{-(1+z)t}\mathrm{d}t\\
 &=& \frac{1}{1+z}e^{-\rho(x)-z(x+\sigma_\eps)},
\eean
with the estimate
\bean
\|\hat{\psi}(z)\|_{X_1} &\le& \frac{1}{|1+z|}\int_0^\infty e^{-(x-\sigma_\eps)_+-\Re e \,z(x+\sigma_\eps)}\mathrm{d}x\\
 &\le& \frac{e^{\sigma_\eps(1-\Re e \,z)}}{|1+z|(1+\Re e \,z)}.
\eean
All together, we get
$$
\|(R_{\BB^1_\eps}(z)\AA_\eps)^2(z)(g,v)\|_{X_1}\le\|\hat{\psi}(x)\|_{X_1}|\RR_\eps[g,v]|\le\frac{C}{\langle z\rangle}\|(g,v)\|_X,
$$
for any $z\in\Delta_{-1}$ and some constant $C>0$ and similar estimate for the second part
$$
\|(R_{\BB^2_\eps}(z)\AA_\eps)^2(z)(g,v)\|_{X_2}\le\frac{C}{\langle z\rangle}\|(g,v)\|_X.
$$
\endproof

\smallskip
\proof[Proof of Theorem~\ref{th:MReFR}]
From the factorization formula with the properties of the two auxiliarty operators in the above Lemma \ref{lem:LmFR} and benefiting from the perturbation argument in \cite{MMcmp, T*, Mbook*}, we deduce that for any $\langle g_0\rangle=0$, there holds
$$
\|g(t)\|_{X_1}=\|S_{\LLL^1_\eps}(t)g_0\|_{X_1}\leq C e^{-t}\,\|g_0\|_{X_1}.
$$
The Duhamel formula associated to the equation $\partial_t v=\LLL^2_\eps (g,v)$ writes
$$
v(t)=S_{\BB^2_\eps}(t) v_0+\int_0^t S_{\BB^2_\eps}(t-s)\AA^2_\eps\big(g(s), v(s)\big)\,\mathrm{d}s.
$$
Using the already known  estimate on $g(t)$, we deduce
\bean
\|S_{\LLL^2_\eps}(t)v_0\|_{X_2} &=& \|v(t)\|_{X_2}\leq \|S_{\BB^2_\eps}(t)
v_0\|_{X_2}+\int_0^t \|S_{\BB^2_\eps}(t-s)\delta_0\NN_\eps[g(s)]\|_{X_2} \,
\mathrm{d}s\\
&\leq& e^{-\delta t} \|v_0\|_{X_2}+\int_0^t e^{-\delta(t-s)} C\,
e^{-s}\|g_0\|_{X_1}\, \mathrm{d}s\\
&\leq& C\, e^{\alpha t}\|(g_0,\,v_0)\|_{X}
\eean
for some $0>\alpha>\max\{-1, -\delta\}$, which yields our
conclusion.
\endproof

\smallskip

Now, we complete the rest part of the proof of Theorem \ref{th:MR} to describe the stability
in the case with delay more precisely compared to that in \cite{PK1}.
\proof[Proof of Theorem~\ref{th:MR} for the step function firing rate in the case with delay]
We write the system as
\bean
\partial_t f &=& - \partial_x f - f\mathbf{1}_{x>\sigma(\eps\DD[u])} + \delta_0 \PPP[f, \eps\DD[u]]\\
\partial_t u &=& - \partial_y u  + \delta_0 \PPP[f,\eps\DD[u]]
\eean
with 
$$
\PPP[f,m]= \int a(m) f, \quad \DD[u] = \int b u.
$$
We recall that the steady state $(F,U)$, where $U:= M {\bf 1}_{y \ge 0}$, satisfies
\bean
0 &=& - \partial_x F- F\mathbf{1}_{x>\sigma(\eps M)} + \delta_0 M  
 \\
 0&=& - \partial_y U  + \delta_0 M, \quad M = \DD[U] = \PPP[F,\eps\DD[U]].
\eean
We introduce the variation $g := f - F$ and $v = u - U$. The equation on $g$ is 
\bean
\partial_t g 
&=& - \partial_x g - f\mathbf{1}_{x>\sigma(\eps\DD[u])}
+F\mathbf{1}_{x>\sigma(\eps M)} + \delta_0 (\PPP[f,\eps\DD[u]] - \PPP[F,\eps\DD[U]])\\
&=& \LLL_\eps^1(g,v) - \QQ[g,v] + \delta_0\langle\QQ[g,v]\rangle\\ 
&=& \LLL_\eps^1(g,v) + \ZZ^1[g,v],
\eean
with
\bean
\ZZ^1[g,v] &=& - \QQ[g,v] + \delta_0\langle\QQ[g,v]\rangle,\\
\QQ[g,v] &:=& \sgn\big(\DD[v]\big)(g+F_\eps)\mathbf{1}_{\III[\DD[v]]},
\eean
where the interval 
$$
\III(n):=(\sigma(\eps M+\eps n_+),\sigma(\eps M+\eps n_-)].
$$
The equation on $v$ is 
\bean
\partial_t v &=& - \partial_y v  + \delta_0 (\PPP[f,\eps\DD[u]] - \PPP[F,\eps\DD[U]])\\ 
&=& - \partial_y v +  \delta_0 \OO[g,v] + \delta_0\langle\QQ[g,v]\rangle\\ 
&=& \LLL_\eps^2(g,v) + \ZZ^2[g,v],
\eean
with
$$
\ZZ^2[g,v] := \delta_0\langle\QQ[g,v]\rangle.
$$
We observe that
\bean
\|\QQ[g,v]\|_{L^p} &=& \|\sgn(\DD[v])g\mathbf{1}_{\III[\DD[v]]}\|_{L^p}\leq\big|\III[\DD[v]]\big|^{1/p}\\
 &=& \big(\sigma(\eps M-\eps\DD[v]_-)-\sigma(\eps M+\eps\DD[v]_+)\big)^{1/p}\\
 &\le& C\big(\eps\|\sigma'\|_\infty\big|\DD[v]\big|\big)^{1/p}\\
 &\le& C\,\eps(\zeta_\eps)\|v\|_{X_2},
\eean
which implies immediately that
\bean
\|\ZZ^1[g,v]\|_{X_1} &\le& \eps\,(\zeta_\eps)\,C\|(g,v)\|_X,\\
\|\ZZ^2[g,v]\|_{X_2} &\le& \eps\,(\zeta_\eps)\,C\|(g,v)\|_X.
\eean
We write the Duhamel formula
$$
\big(g(t),v(t)\big) = S_{\LLL_\eps}(t) (g_0,v_0) + \int_0^t S_{\LLL_\eps}(t-s) \ZZ[g(s),v(s)] \, \mathrm{d}s.
$$
and thanks to the Gronwall's Lemma, we conclude the exponential asymptotic stability of Theorem \ref{th:MR} for the step function firing rate in the case with delay.
\endproof

\bigskip
\bibliographystyle{acm}

\signqw  

\end{document}